\documentclass[12pt]{article}

\usepackage{url,epsf,amsmath,amssymb,amsfonts}
\usepackage{color}

\newtheorem{theorem}{\bf Theorem}
\newtheorem{corollary}[theorem]{\bf Corollary}

\newcommand{\cp}{\,\square\,}
\newcommand{\gp}{{\rm gp}}
\newcommand{\pn}{P_\infty^n}

\begin{document}

\title{The general position number of integer lattices}

\author{
Sandi Klav\v zar $^{a,b,c}$
\and
Gregor Rus $^{c,d}$
}

\date{}

\maketitle

\begin{center}
$^a$ FMF, 
University of Ljubljana, Slovenia \\
\medskip

$^b$ FNM, 
University of Maribor, Slovenia \\
\medskip

$^{c}$ IMFM, 
Ljubljana, Slovenia \\
\medskip

$^d$ FOV, 
University of Maribor, Slovenia\\

\end{center}

\begin{abstract}
The general position number ${\rm gp}(G)$ of a connected graph $G$ is the cardinality of a largest set $S$ of vertices such that no three pairwise distinct vertices from $S$ lie on a common geodesic. The $n$-dimensional grid graph $\pn$ is the Cartesian product of $n$ copies of the two-way infinite path $P_\infty$. It is proved that if $n\in {\mathbb N}$, then ${\rm gp}({P_\infty^n}) = 2^{2^{n-1}}$. The result was earlier known only for $n\in \{1,2\}$ and partially for  $n=3$.  
\end{abstract}

\noindent
{\bf E-mails}: sandi.klavzar@fmf.uni-lj.si, gregor.rus4@um.si


\medskip\noindent
{\bf AMS Subj. Class.}: 05C12, 05C76, 11B75

\section{Introduction and preliminaries}
\label{sec:intro}

Subsets of vertices of (infinite) grids with special properties are of wide interest, the variety~\cite{bouznif-2019, distefano-2017, drews-2019, dyb-2020, gan-2018, guo-2020, jobsen-2018}  of such investigations clearly supports this statement. In this note we are interested in largest general position sets of grids. By now only some specific results about these sets in grids were known, here we solve the problem completely.

A set $S$ of vertices of a connected graph $G$ is a {\em general position set} if $d_G(u,v) \ne d_G(u,w) + d_G(w,v)$ holds for every $\{u,v,w\}\in \binom{S}{3}$, where $d_G(x,y)$ denotes the shortest-path distance between $x$ and $y$ in $G$. The {\em general position number} $\gp(G)$ of $G$ is the cardinality of a largest general position set in $G$.  This concept and terminology was introduced in~\cite{manuel-2018a}, in part motivated by the century old Dudeney's No-three-in-line problem~\cite{dudeney-1917}. A couple of years earlier and in different terminology, the problem was also considered in~\cite{ullas-2016}. Moreover, in the special case of hypercubes, the general position problem has been studied back in 1995 by K\"orner~\cite{korner-1995}. Following these seminal papers, the general position problem has been studied from different perspectives in several subsequent papers~\cite{anand-2019, ghorbani-2019, klavzar-2019+, klavzar-2019b, manuel-2018b, patkos-2019+}. 

The {\em Cartesian product} $X\cp Y$ of graphs $X$ and $Y$ is the graph with the vertex set $V(X) \times V(Y)$, vertices $(x,y)$ and $(x',y')$ being adjacent if either $x=x'$ and $yy'\in E(Y)$, or $y=y'$ and $xx'\in E(X)$. The  Cartesian product $X_1\cp \cdots \cp X_n$, where each factor $X_i$ is isomorphic to $X$, will be shortly denoted by $X^n$. If $P_\infty$ denotes the two-way infinite path, then one of the main results from~\cite{manuel-2018b} asserts that $\gp(P_\infty^2) = 4$. In the same paper it was also proved that $10\le \gp(P_\infty^3) \le 16$. The lower bound $10$ was improved to $14$ in~\cite{klavzar-2019b}. In this note we round these investigations by proving the following result. 

\begin{theorem}
\label{thm:main}
If $n\in {\mathbb N}$, then $\gp(\pn) = 2^{2^{n-1}}$.
\end{theorem}

In the rest of the section we list some further notation and and preliminary results. In the next section we prove the theorem. In the concluding section we give a couple of consequences of the theorem and pose an open problem.  

For a positive integer $k$ we will use the notation $[k] = \{1,\ldots, k\}$. Throughout we will set $V(P_\infty) = {\mathbb Z}$, where $u,v\in V(P_\infty)$ are adjacent if and only if $|u-v| = 1$. With this convention we have $V(P_\infty^n) = {\mathbb Z}^n$. If $u\in V(\pn)$, then for the coordinates of $u$ we will use the notation $u = (u_1, \ldots, u_n)$. If a vertex from  $V(\pn)$ will be indexed, say $u_i\in V(\pn)$, then this notation will be extended as $u_i = (u_{i,1}, \ldots, u_{i,n})$. From the Distance Lemma~\cite[Lemma 12.2]{imrich-2008} it follows that
\begin{equation}
\label{eq:distance}
d_{\pn}(u, v) = \sum_{i=1}^n |u_i - v_i|\,. 
\end{equation}  
From here it is not difficult to deduce that a vertex $w\in V(\pn)$ lies on a shortest $u,v$-path in $\pn$ if and only if  $\min \{u_i,v_i\}\leq w_i\leq \max\{u_i,v_i\}$ holds for every $i\in [n]$.  

A sequence of real numbers is {\it monotone} if it is monotonically increasing or monotonically decreasing. The celebrated Erd\H{o}s-Szekeres result on monotone sequences reads as follows (cf.\ also~\cite[Theorem 1.1]{bukh-2014}). 

\begin{theorem} {\rm \cite{ErSz35}}
\label{thm:Er-Sze}
For every $n\ge 2$, every sequence $(a_1,\ldots, a_N)$ of real numbers with $N\ge (n-1)^2 + 1$ elements contains a monotone subsequence of length $n$. 
\end{theorem}

\section{Proof of Theorem~\ref{thm:main}}
\label{sec:proof}

Theorem is obviously true for $n = 1$ and was proved for $n=2$ in~\cite[Corollary 3.2]{manuel-2018b}. 

Let now $n\ge 3$ and let $U^{(1)} = \{u_1, \ldots, u_{2^{2^{n-1}} +1}\}$  be a set of vertices of $\pn$ of cardinality $2^{2^{n-1}} +1$. We may without loss of generality assume that the first coordinates of the vertices from $U^{(1)}$ are ordered, that is, $u_{1,1}\le u_{2,1} \le \cdots \le u_{{2^{2^{n-1}} +1},1}$. By Theorem~\ref{thm:Er-Sze}, there exists a subset $U^{(2)}$ of $U^{(1)}$ of cardinality $2^{2^{n-2}} +1$, such that the second coordinates of the vertices from $U^{(2)}$ form a monotone (sub)sequence. Inductively applying this argument, we arrive at a set $U^{(n)}\subset U^{(n-1)}$ of cardinality $2^{2^{n-n}} +1 = 3$, in which the $n^{\rm th}$ coordinates of the three vertices form a monotone (sub)sequence. As $U^{(n)}\subset U^{(n-1)} \subset \cdots \subset U^{(1)}$, the induction argument yields that for every $i\in [n-1]$, the $i^{\rm th}$ coordinates of the vertices from $U^{(n)}$ likewise form a monotone (sub)sequence. If $U^{(n)} = \{u,v,w\}$, where $u_1 \le v_1 \le w_1$, this implies (having \eqref{eq:distance} in mind) that $v$ lies on a shortest $u,w$-path. We conclude that $\gp(\pn) \le 2^{2^{n-1}}$.

To prove the other inequality we are going to inductively construct a general position set $X^{(n)} = \{x_1^{(n)}, \ldots, x_{2^{2^{n-1}}}^{(n)} \}$ for $n\ge 2$ as follows. Set $X^{(2)} = \{(1,2),(2,1),(3,4),(4,3)\}$, where $x_1^{(2)} = (1,2)$, $x_2^{(2)} = (2,1)$, $x_3^{(2)} = (3,4)$, and $x_4^{(2)} = (4,3)$.  Suppose now that $X^{(n-1)}$ is defined for some $n\ge 3$, and construct $X^{(n)}$ as follows. Set first 
\begin{itemize}
\item $x_{i,1}^{(n)} = i$, $i\in [2^{2^{n-1}}]$.
\end{itemize}
Next, write each $i\in [2^{2^{n-1}}]$ as $i = p\cdot 2^{2^{n-2}} + r$, where $0\le p < 2^{2^{n-2}}$ and $r\in [2^{2^{n-2}}]$, and set  
\begin{itemize}
\item $x_{i,j}^{(n)} = \left(x_{(p+1),j}^{(n-1)} - 1\right)\cdot 2^{2^{n-2}} + x_{r,j}^{(n-1)}$ for each $j\in \{2,\ldots, n-1\}$, and 
\item $x_{i,n}^{(n)} = x^{(n-1)}_{p+1,n-1}\cdot 2^{2^{n-2}} - x^{(n-1)}_{r,n-1}+1\,.$
\end{itemize}
Roughly speaking, for the $j^{\rm th}$ coordinate, where $j\in \{2,\ldots, n-1\}$, we partition the sequence $(x^{(n)}_{i,j})_{i=1}^{2^{2^{n-1}}}$ into $2^{2^{n-2}}$ blocks each of $2^{2^{n-2}}$ values and sort the blocks as well as the values inside the blocks according to the values $(x^{(n-1)}_{i,j})_{i=1}^{2^{2^{n-2}}}$. The values of the $n^{\rm th}$ coordinate is then obtained from the values of the $(n-1)^{\rm th}$ coordinate by reversing the sequence in each of the $2^{2^{n-2}}$ blocks, while keeping the sequence of the blocks. For example, the coordinates of the vertices from $X^{(3)}$ are shown in Table~\ref{tab:coordinates}. 

 \begin{table}[ht]
 \centering      
\begin{tabular}{|c||c|c|}
\hline
$i$ & 1 & 2\\
\hline \hline
$x^{(1)}_{i,1}$ & 1&2 \\
\hline
\end{tabular}

\bigskip

\begin{tabular}{|c||c|c|c|c|}
\hline
$i$ & 1 & 2 & 3 & 4\\
\hline \hline
$x^{(2)}_{i,1}$ & 1&2 & 3 & 4 \\
\hline
$x^{(2)}_{i,2}$ & 2 & 1 & 4 & 3 \\
\hline
\end{tabular} 

\bigskip

      \begin{tabular}{|c||c|c|c|c|c|c|c|c|c|c|c|c|c|c|c|c|c|}
\hline
$i$ & 1 & 2 & 3 & 4 & 5 & 6 & 7 & 8 & 9 & 10 & 11 & 12 & 13 & 14 & 15 & 16 \\
\hline\hline
$x^{(3)}_{i,1}$ & 1 & 2 & 3 & 4 & 5 & 6 & 7 & 8 & 9 & 10 & 11 & 12 & 13 & 14 & 15 & 16 \\
\hline
$x^{(3)}_{i,2}$ & 6 & 5 & 8 & 7 & 2 & 1 & 4 & 3 & 14 & 13 & 16 & 15 & 10 & 9 & 12 & 11 \\
\hline
$x^{(3)}_{i,3}$ & 7 & 8 & 5 & 6 & 3 & 4 & 1 & 2 & 15 & 16 & 13 & 14 & 11 & 12 & 9 & 10 \\
\hline
 \end{tabular}
 

     \caption{Coordinates of the vertices in sets $X^{(1)}$, $X^{(2)}$, and $X^{(3)}$, respectively.}
     \label{tab:coordinates}
 \end{table}

To complete the proof it suffices to show that for each $n\ge 2$, the set $X^{(n)}$ forms a general position set of $\pn$. We proceed by induction on $n$, the base case $n=2$ being clear. Suppose now that $X^{(n-1)}$ is a general position set of $P_\infty^{n-1}$ and consider $X^{(n)}$. Partition the set $[2^{2^{n-1}}]$ into $2^{2^{n-2}}$ {\em blocks} $\{1,\ldots, 2^{2^{n-2}}\}$, $\{2^{2^{n-2}}+1,\ldots, 2^{2^{n-2}+1}\}$, $\ldots$  Let $x_i^{(n)}$, $x_j^{(n)}$, and $x_k^{(n)}$ be pairwise different vertices from  $X^{(n)}$ and consider the following three cases, where again write each number $m\in [2^{2^{n-1}}]$ as $m = p_m\cdot 2^{2^{n-2}} + r_m$, where  $0\le p_m < 2^{2^{n-2}}$ and $r_m\in [2^{2^{n-2}}]$.

To prove that no three vertices of $X^{(n)}$ lie on a common geodesic, the following claim will be useful. 

\medskip\noindent
{\bf Claim A} If $n\ge 3$ and if $i$ and $j$ are in the same block, then  $x^{(n)}_{i,n-1}<x^{(n)}_{j,n-1}$ if and only if $x^{(n)}_{i,n}>x^{(n)}_{j,n}$.

\medskip\noindent
Let $n\ge 3$ and let $i = p_i \cdot 2^{2^{n-2}} + r_i$ and $j = p_j \cdot 2^{2^{n-2}} + r_j$. Assume  that  $x^{(n)}_{i,n-1}<x^{(n)}_{j,n-1}$. By the construction, $x^{(n)}_{i,n-1} = \left(x_{(p_i+1),n-1}^{(n-1)} - 1\right)\cdot 2^{2^{n-2}} + x_{r_i,n-1}^{(n-1)}$ and $x^{(n)}_{j,n-1} = \left(x_{(p_j+1),n-1}^{(n-1)} - 1\right)\cdot 2^{2^{n-2}} + x_{r_j,n-1}^{(n-1)}$. Since $i$ and $j$ are in the same block, that is, $p_i=p_j$, we get that  $x^{(n)}_{j,n-1}-x^{(n)}_{i,n-1} = x_{r_j,n-1}^{(n-1)} - x_{r_i,n-1}^{(n-1)}$. As we have assume that $x^{(n)}_{i,n-1}<x^{(n)}_{j,n-1}$, it follows that $x_{r_j,n-1}^{(n-1)} > x_{r_i,n-1}^{(n-1)}$. Since $x^{(n)}_{i,n} = \left(x_{(p_i+1),n-1}^{(n-1)} \right)\cdot 2^{2^{n-2}} - x_{r_i,n-1}^{(n-1)} + 1$ and $x^{(n)}_{j,n} = \left(x_{(p_j+1),n-1}^{(n-1)} \right)\cdot 2^{2^{n-2}} - x_{r_j,n-1}^{(n-1)} + 1$, we conclude that $x^{(n)}_{i,n}-x^{(n)}_{j,n} = x_{r_j,n-1}^{(n-1)} - x_{r_i,n-1}^{(n-1)} > 0$. The reverse implication is proved along the same lines. This proves Claim~A. 

\medskip
We now distinguish the following cases.  

\medskip\noindent
{\bf Case 1}: $|\{p_i, p_j, p_k\}| = 1$. \\
In this case $i$, $j$, and $k$ belong to the same block. Then by the induction hypothesis, the first $n-1$ coordinates assure that $x_i^{(n)}$, $x_j^{(n)}$, and $x_k^{(n)}$ do not lie on a common geodesic.

\medskip\noindent
{\bf Case 2}: $|\{p_i, p_j, p_k\}| = 2$.\\
In this case we may assume without loss of generality that $p_i = p_j  < p_k$. Further assuming without loss of generality that  $r_i < r_j$, we have $x^{(n)}_{i,1} < x^{(n)}_{j,1} < x^{(n)}_{k,1}$. Hence if $x_i^{(n)}$, $x_j^{(n)}$, and $x_k^{(n)}$  lie on a common geodesic, then necessarily $x_j^{(n)}$ lies between $x_i^{(n)}$ and $x_k^{(n)}$. Suppose first that $x^{(n)}_{i,n-1} < x^{(n)}_{j,n-1} < x^{(n)}_{k,n-1}$. Then by Claim~A we have $x^{(n)}_{j,n} < x^{(n)}_{i,n}$. Since $x^{(n)}_{i,n-1} < x^{(n)}_{k,n-1}$ and $i, k$ are in different blocks, we have $x^{(n)}_{i,n} < x^{(n)}_{k,n}$. Therefore, $x^{(n)}_{j,n} < x^{(n)}_{i,n} < x^{(n)}_{k,n}$, so $x_i^{(n)}$, $x_j^{(n)}$, and $x_k^{(n)}$ do not lie on a common geodesic. By a parallel argument we come to the same conclusion if $x^{(n)}_{i,n-1} > x^{(n)}_{j,n-1} > x^{(n)}_{k,n-1}$ holds. 

\medskip\noindent
{\bf Case 3}: $|\{p_i, p_j, p_k\}| = 3$.\\
Since $p_i$, $p_j$ and $p_k$ are pairwise different, which means then we can consider their blocks. Considering the whole blocks as single contracted vertices, their components are the integral part of dividing components with $2^{2^{n-2}}$ (the $p$ values of the components). Since $p$ values are obtained in the construction from a general position set of the $(n-1)$ dimensional grid, the induction hypothesis assures that these contracted vertices do not lie on a common geodesic, since this contracted vertices are part of a general position set in $P_\infty^{n-1}$. This in turn implies that also $x_i^{(n)}$, $x_j^{(n)}$, and $x_k^{(n)}$ do not lie on a common geodesic.

\section{Concluding remarks}
\label{sec:conclude}

Recall that a subgraph $H$ of a graph $G$ is {\em isometric} if $d_H(u,v) = d_G(u,v)$ holds for each pair of vertices $u,v\in V(H)$. Since $P_{i_1}\cp \cdots \cp P_{i_n}$ is an isometric subgraph of $\pn$, Theorem~\ref{thm:main} immediately implies:

\begin{corollary} 
\label{cor:finite-grids}
If $n\ge 2$, and $i_1, \ldots, i_n\ge 2^{2^{n-1}}$, then $\gp(P_{i_1}\cp \cdots \cp P_{i_n}) = 2^{2^{n-1}}$. 
\end{corollary} 

More generally, if a graph $G$ contains an isometric grid $P_{i_1}\cp \cdots \cp P_{i_n}$, where each $i_j\ge 2^{2^{n-1}}$, then $\gp(G) \ge 2^{2^{n-1}}$. For instance: 

\begin{corollary} 
\label{cor:finite-turus}
If $n\ge 2$, and $i_1, \ldots, i_n\ge 2^{2^{n-1}+1}$, then $\gp(C_{i_1}\cp \cdots \cp C_{i_n}) \ge 2^{2^{n-1}}$. 
\end{corollary} 

From~\cite{ghorbani-2019} we know that $\gp(G\cp H) \geq \gp(G) + \gp(H) - 2$ holds for finite, connected graphs $G$ and $H$.  Since the general position number of a path is $2$, Corollary~\ref{cor:finite-grids} demonstrates  that the difference in the inequality can be arbitrary large. 

In~\cite{klavzar-2019b} a formula for the number of general position sets of cardinality $4$ in $P_r\cp P_s$ (that is, of largest general position sets) is determined for each $r, s\ge 2$. Because of this result and Corollary~\ref{cor:finite-grids}, an interesting and intriguing problem is to determine the number of largest general position sets in $P_{i_1}\cp \cdots \cp P_{i_n}$, where $n\ge 3$ and $i_1, \ldots, i_n\ge 2^{2^{n-1}}$.

\section*{Acknowledgements}

We acknowledge the financial support from the Slovenian Research Agency (research core funding No.\ P1-0297 and projects J1-9109, J1-1693, N1-0095, N1-0108).

\end{document}